\documentclass{amsproc}
\usepackage{bm,color}

\usepackage{graphicx}

\newtheorem{theorem}{Theorem}[section]
\newtheorem{lemma}[theorem]{Lemma}

\theoremstyle{definition}
\newtheorem{definition}[theorem]{Definition}

\newtheorem{corollary}[theorem]{Corollary}
\newtheorem{proposition}[theorem]{Proposition}
\theoremstyle{remark}
\newtheorem{remark}[theorem]{Remark}

\numberwithin{equation}{section}

\def\bmtx{\begin{matrix}}
\def\emtx{\end{matrix}}

\def\ovsig{\overline{\sigma}}

\def\wb{[\bfb]}

\def\bfb{{\mathbf b}}

\def\bfu{{\mathbf u}}
\def\bfv{{\mathbf v}}
\def\bfw{{\mathbf w}}

\def\d{\partial}

\def\sgn{\mathrm{sgn}}

\def\ZZ{\mathbb Z}

\def\RR{\mathbb R}

\def\QQ{\mathbb Q}
\def\1{{\mathbf 1}}
\def\0{{\mathbf 0}}

\def\cocoa{{\hbox{\rm C\kern-.13em o\kern-.07em C\kern-.13em o\kern-.15em A}}}

\def\Dcal{\mathcal D}

\def\Span{{\rm Span}}

\def\ttp{{\tt p}}

\def\xx{{\bf x}}

\def\bfu{{\bf u}}

\def\ovDc{\overline{\Dcal}}
\def\ovD{\overline{D}}

\def\End{\mathrm{End}}

\def\blamb{{\bm \lambda}}

\def\Pcal{{\mathcal P}}

\def\w2M{\bigwedge^2M}

\def\wM{\bigwedge M}

\def\w{\wedge }
\def\bw{\bigwedge }

\protect
\protect
\protect
\protect

\protect
\protect

\def\sra{\rightarrow}

\def\proof{\noindent{\bf Proof.}\,\,}
\def\qed{{\hfill\vrule height4pt width4pt depth0pt}\medskip}
\def\be{\begin{equation}}
\def\ee{\end{equation}}

\def\bclm{\begin{claim}}
\def\eclm{\end{claim}}
\def\beqn{\begin{eqnarray}}
\def\eeqn{\end{eqnarray}}
\def\beqn*{\begin{eqnarray*}}
\def\eeqn*{\end{eqnarray*}}

\title[Hasse--Schmidt Derivations and Cayley--Hamilton Theorem ] {Hasse--Schmidt Derivations and Cayley--Hamilton Theorem for Exterior Algebras}

\author[Gatto]{Letterio Gatto}
\address{Dipartimento di Scienze Matematiche\\
Politecnico di Torino\\
C.so Duca degli Abruzzi, 24\\
10129 - Torino, Italia}
\curraddr{}
\email{letterio.gatto@polito.it}
\thanks{Work parially sponsored by PRIN ``Geometria sulle Variet\`a algebriche'', INDAM-GNSAGA e "Finanziamento diffuso della Ricerca" del Politecnico di Torino. }

\author[Scherbak]{Inna Scherbak}
\address{School of Mathematical Sciences\\
Raymond \& Beverly Sackler 
Faculty of Exact Sciences\\
Tel Aviv University, P.O. Box 39040, Tel Aviv 6997801, Israel}
\curraddr{}
\email{scherbak@post.tau.ac.il}
\thanks{}

\subjclass[2010]{Primary 15A75, 17B69. 050E5}

\dedicatory{In memory of those participants of the Voronezh Winter Mathematical School\\
who have already passed into another world where all problems are solved}

\keywords{Hasse-Schmidt Derivations on Grassmann Algebras, Theorem of Cayley and Hamilton, Vertex Operators}

\date{}
\begin{document}

\begin{abstract}
Using the natural notion of  {\em Hasse--Schmidt derivations on an exterior algebra}, we relate two classical and seemingly unrelated subjects. The first  is the famous Cayley--Hamilton theorem of linear algebra,
``{\em each endomorphism of a finite-dimensional vector space is a root of its own characteristic polynomial}'',  and the second concerns the expression of the bosonic vertex operators occurring in the representation theory of the (infinite-dimensional) Heinsenberg algebra.  
\end{abstract}
\maketitle

\section {Introduction} \label{S1}
In 1937, Hasse and Schmidt  introduced the notion of higher derivations \cite{HaS}, nowadays called {\sl Hasse--Schmidt (HS) derivations}. 
Let $(A,*)$ be an algebra over a ring $B$, not necessarily commutative or associative. A {\sl HS--derivation on} $A$ is a $B$-algebra homomorphism, $D(t): A\sra A[[t]]$, that is, a $B$-linear mapping satisfying
$$ D(t)(a_1*a_2)=D(t)a_1*D(t)a_2,\ \ \forall a_1,a_2\in A.$$

A  fundamental example of a HS--derivation is given by the map 
sending any function $f=f(z)$, holomorphic in  some domain of the complex plane, to its {\sl formal Taylor series}, 
$$
f(z)\mapsto T(t)[f(z)]=\left(\exp\left(t{d\over dz}\right)\right)f(z). 
$$
The property $T(t)[f(z)g(z)]=T(t)[f(z)]T(t)[g(z)]$ encodes  the full set of the Leibniz's rules, 
$$
{d^i(fg)\over dz^i}=\sum_{j=0}
^i\,{i\choose j}{d^jf\over dz^i}\cdot {d^{i-j}f\over dz^{i-j}}\,.
$$

In general,  if $(A,*)$ is any commutative $\QQ$-algebra and  $\delta(t):A\sra A[[t]]$ is a derivation in the  Leibniz rule sense, i.e., $\delta(t)(a*b)=\delta(t)a*b+a*\delta(t)b$, then $\exp (\delta(t))$ is a HS--derivation. 

The aim of Hasse and Schmidt was to find a counterpart to the Taylor series that would work in positive characteristic. Their definition does not require division by integers and is therefore particularly suitable for this purpose. Schmidt later applied  the theory to investigate Weierstrass points and Wronskians  on curves in positive characteristic~ \cite{Sch}.

\medskip
In a number of papers motivated by Schubert Calculus \cite{G1,SCGA} (see also the book  \cite{GatSal3}),  one of us proposed to study  HS--derivations for exterior algebras.

If $A$ is a commutative ring with unit, $M$ a module over $A$,  and $\bw M$  its exterior algebra, then a {\sl HS--derivation on} $\bw M$ is a $\wedge$-homomorphism 
$$\Dcal(t): \bw M\sra \left(\bw M\right)[[t]],$$
that is, a linear mapping satisfying
\be
\Dcal(t)(u\w v)=\Dcal(t)u\w \Dcal(t)v,\ \ \forall  u,v\in\bw M. \label{eq:Duv}
\ee

\medskip
In this paper we consider HS--derivations  $\Dcal(t)=\sum_{i\geq 0} \Dcal_i\cdot t^i$ 
with $\Dcal_0=\1$, the identity on $\bw M$.  Such $\Dcal(t)$ is  invertible as an element of $\End{(\bw M)}[[t]]$, that is, there exists $\ovDc(t)\in\End{(\bw M)}[[t]]$ satisfying
\be
\ovDc(t)D(t)=\Dcal(t)\ovDc(t)=\1. \label{eq1:eq1}
\ee
A straightforward calculation shows that $\ovDc(t)$ is also a HS--derivation on $\bw M$. Hence,
\be
\Dcal(t)(\ovDc(t)u\w v)=u\w \Dcal(t)v,\  \ovDc(t)(\Dcal(t)u\w v)=u\w \ovDc(t)v,\ \ \forall u,v\in\bw M .
\label{eq:1}
\ee
The coefficients  of $ t $ in these equations give
$u\w \Dcal_1v=\Dcal_1(u\w v)-\Dcal_1u\w v$. That  
 is why in~\cite{G1} we call~(\ref{eq:1}) the {\sl integration by parts} formulas.

\bigskip
In the present article we show how these simple formulas link  two classical and seemingly unrelated subjects (one  finite-dimensional and the other infinite-dimensional), apparently leading to a unified  interpretation.
 
 \medskip
One topic is the classical  Cayley--Hamilton Theorem of Linear Algebra saying that  {\sl each endomorphism $f$ of an $r$-dimensional vector space $M$ is a root of its own characteristic polynomial $\det(t\1-f)$}. Let us reformulate this theorem  as a linear recurrence relation on the sequence of endomorphisms $(f^j)_{j\geq 0}$,
\be f^{r+k}-e_1f^{r+k-1}+\cdots+(-1)^ke_rf^k=\0,\ \ \forall k\geq 0,\label{eqn:clCH}
\ee
where $\det(t\1-f)=t^r-e_1t^{r-1}+\cdots+(-1)^re_r$, and $\0$ denotes the zero endomorphism.
Consider $D(t)=\sum_{i\geq 0} D_i\cdot t^i$,  the unique HS--derivation on the exterior algebra of $M$ such that ${D_i}_{|M}=f^i$, $i\geq 0$. It turns out that the sequence $(D_i)_{i\geq 0}$ of endomorphisms of $\bw M$ satisfies  relations similar to~(\ref{eqn:clCH}), see  Theorem~\ref{thm23} in Section~\ref{S21} for the exact formulation, and Section~\ref{S3} for the proof.

\medskip
The other topic concerns bosonic vertex operators arising in the representation theory of the (infinite-dimensional) Heisenberg algebra (see, for example, \cite{Ka}). As we observe in Section~\ref{S22}, any  countably generated vector space over the rationals can be equipped with the structure of a  free module of  finite rank $r$ over a ring of polynomials in $r$ variables  with rational coefficients, for any integer $r>0$. We present the construction in Section~\ref{S4}, and in Section~\ref{S5} we apply it to obtain the ``finite-dimensional approximation''  to the well-known expressions of the vertex operators $\Gamma(t)$ and $\Gamma^*(t)$  generating  the bosonic Heisenberg vertex algebra (see \cite[p.~56]{Ka}). We interpret $\Gamma(t)$ and $\Gamma^*(t)$  as the limit, when $r\sra\infty$, of the ratio of two characteristic polynomials associated to the shift endomorphisms of steps $+1$ and $-1$, respectively. The precise  formulation can be found in Section~\ref{S23}.

\medskip
 Our work is based on interpreting~ (\ref{eq:1}) as a sort of abstract Cayley--Hamilton theorem, holding for general invertible HS--derivations on exterior algebras of  arbitrary modules (not necessarily free).  If the module is countably generated, then~(\ref{eq:1}) produces a sequence of Cayley-Hamilton relations~(\ref{f25}) which specialize to the classical Cayley--Hamilton formulas~(\ref{eqt22:prech}) when considering Hasse--Schmidt derivations associated to endomorphisms of finitely generated free modules.

\medskip\noindent
{\bf Plan of the paper.}   In Section~\ref{S2} we formulate the  main statements.   Section~\ref{S21} is devoted to  our  extention of the Cayley--Hamilton theorem on an exterior algebra of  a finite rank free module, which is  Theorem~\ref{thm23}.  The proof can be found in Section~\ref{S3}, which also includes  the necessary information on the Hasse--Schmidt derivations,   and the discussion concerning the case when the ring $A$ contains the rationals.

In Section~\ref{S22},  we equip  a countably infinite-dimensional  $\QQ$-vector  space  with a natural structure of a  free module of rank $r$ over the ring of polynomials of $r$ variables  with rational coefficients, for any integer $r>0$. The construction is based on a Giambelli's type formula. The details are explained in Section~\ref{S4}. 

In Section~\ref{S23} we apply our construction to the bosonic Heisenberg vertex algebra and interpret the truncation of bosonic vertex operators as the ratio of two characteristic polynomials,  respectively associated to the shift endomorphisms of step $\pm 1$; see  Section~\ref{S5} for detailed explanation.

\medskip\noindent
{\bf Acknowledgments} We thank the referees  for their efforts reading the first version of the paper. Their sharp criticism really contributed to the improvement of the article. In particular, we appreciate the help of the referee who found a mistake in our text. Her/his suggestions allowed us to simplify our proof. We are also grateful to the other referee who demanded  more motivations.

\section{Formulation of the results}\label{S2}
\subsection{Cayley--Hamilton theorem for exterior algebras}\label{S21}
Let $M$ be a free $A$-module of at most countable (i.e., either finite or countable) rank.
If $b_1, b_2,\ldots $ is a basis of $M$,
then $\left\{b_{i_1}\w\ldots\w b_{i_j}\right\}_{1\leq i_1<\ldots <i_j}$ is a basis of $\bw^jM$, $j\geq 1$. Following a referee's suggestion, we observe that any 
$A$-endomorphism $f:M\sra M$  naturally extends to the endomorphism $\bw f$ of $\bw M$ as follows. 

Define 
$\bw f=\bigoplus_{j\geq0}\bw^j f$ by setting  $\bw^0f$ to be the identity on $\bw^0M=A$, and by defining the  action of $\bw^jf$ on the basis of $\bw^jM$ as follows,
\be 
b_{i_1}\w\ldots\w b_{i_j}\ \mapsto \ fb_{i_1}\w\ldots\w fb_{i_j},\ \  j\geq 1.\label{eqn:basis}
\ee
As a consequence, any formal power series $P(t)=\sum_{i\geq 0}p_it^i$, where $p_i\in\End_AM$,  clearly defines a HS--derivation $\bw P(t)$ on $\bw M$. If, in addition, $p_0=\1$, the identity on $M$, then the formal inverse $P^{-1}(t)\in\left(\End_AM\right)[[t]]$ defines the inverse HS--derivation, 
$$\overline{\bw P}(t)=\bw P^{-1}(t).$$
\begin{proposition}\label{prop11} Let $M$ be an $A$-module of at most countable rank, and $f\in\End_AM$. Set $f^0=\1$. Then
\be  D(t)=\bw\left(\sum_{i\geq 0}f^it^i\right)\label{eqn:D}\ee
is an invertible HS--derivation, and its inverse is
\be \ \ \ \ \ \ \ \ovD(t)=\bw\left(\1-ft\right).\label{eqn:ovD} \qquad \ \ \ \qed \ee
\end{proposition}

\medskip
Our extension of the Cayley--Hamilton theorem concerns modules of finite rank.
If $M$ has rank $r$ over $A$, then $\bw^rM$ has rank 1, and so $D(t)$ and $\ovD(t)$ act on
$\bw^rM$ by multiplication by some formal power series.  Indeed, $\bw^rM=\Span_A\left\{b_1\w\dots\w b_r\right\}$, and, by definition~(\ref{eqn:ovD}),
$$\ovD(t)\left(b_1\w\dots\w b_r\right)=(\1-ft)b_1\w\dots\w(\1-ft)b_r=\det(\1-ft)\left(b_1\w\dots\w 
b_r\right).$$
Thus, the ``eigenvalue'' of $\ovD(t)$ on $\bw^rM$ is $\det(\1-ft)$, considered as an element of $ A[[t]]$. Let us write
$ \det(\1-ft)=E_r(t)$, where
\be E_r(t)=1-e_1t+\cdots+(-1)^rt^r. \label{eqn:E}
\ee
 
\begin{remark}\label{Rem13} 
(1)\ Clearly  ``the eigenvalue''  $H_r(t)$ of $D(t)$ on $\bw^rM$  is the formal inverse of $E_r(t)$, i.e.,
\be  H_r(t)E_r(t)=1.\label{eqn:H} \ee
If we write $H_r(t)=\sum_{j\geq 0}h_j t^j$, then (\ref{eqn:E}) and  (\ref{eqn:H}) determine each $h_j$ as a polynomial of $e_1,\ldots,e_r$. For example, $h_0=1$, $h_1=e_1$, $h_2=e_1^2-e_2$ etc. 

(2)\ It is worth emphasizing that $\{e_i \}$ and $\{h_j \}$ are related in exactly the same way as {\sl the elementary} and {\sl the  complete symmetric functions} of $ r $ variables, since $E_r (t)$ and $H_r(t)$ are their generating functions, respectively (see, for example, \cite [I\,2] {MacDonald}).
\end{remark}

\medskip
For $f\in\End_AM$, one can write HS--derivations $D(t)$ and $\ovD(t)$, defined by~(\ref{eqn:D}) and~(\ref{eqn:ovD}), respectively, in the form
$$D(t)=\sum_{i\geq 0} D_i\cdot t^i,\ \ \ovD(t)=\sum_{i\geq 0}(-1)^i\ovD_i\cdot t^i.$$
 We have $D_0=\1=\ovD_0$,  $D_1=\ovD_1$, and $ D_i\mid_M=f^i$, $i\geq 0$. 

Denote $\bw^{>i}M=\bigoplus_{j=i+1}^r\bw^jM$.
\begin{theorem}\label{thm23}  For\  $1\leq k < r$, the endomorphism 
\be
D_k-e_1D_{k-1}+\cdots+(-1)^ke_k\1 \label{eq2:prech}
\ee
vanishes on $\bw^{>(r-k)}M$, and  for\  $i\geq r$ the endomorphism
 \be
D_i-e_1D_{i-1}+\cdots+(-1)^re_rD_{i-r}\label{eqt22:prech}
 \ee
 vanishes on the whole of $\bw M$.
\end{theorem}

According to~(\ref{eqn:E}), the characteristic polynomial of $f$ is 
\be\det(t\1-f)=t^rE_r(1/t)=t^r-e_1t^{r-1}+\cdots+(-1)^re_r.\label{eqn:1E}\ee
Hence, for $i=r+k$, the restriction of (\ref{eqt22:prech}) to $M$ gives the classical Cayley--Hamilton theorem ~(\ref{eqn:clCH}).

\subsection{A look at the infinite-dimensional case.} \label{S22}
Let $M_0$ be a  $\QQ$-vector space  with a  countable basis, 
and $\bw M_0=\bigoplus_{j\geq 0}\bw^jM_0$ its exterior algebra.

We equip $M_0$ with a structure of a  free module of rank $r$ over the ring of polynomials of $r$ variables  with rational coefficients, for any integer $r>0$.  See Section~\ref{S4} for details.

\medskip
We fix  a basis $(b_j)_{j\geq 1}$ of $M_0$, and  define the shift operators $\sigma_{+1}$, $\sigma_{-1}$ on $M_0$ by their action on the basis,
$$\sigma_{+1}(b_j)= b_{j+1},\ j\geq 1, {\rm \ and\  }\sigma_{-1}(b_1)=0,\ \sigma_{-1}(b_j)=b_{j-1},\ j>1.$$
One can attach to each of the endomorphisms $\sigma_{+1}$, $\sigma_{-1}$  a  unique HS--derivation and its inverse, as in (\ref{eqn:D}), (\ref{eqn:ovD}). In this subsection we need  only the HS--derivations 
$$\sigma_{+}(t), \overline{\sigma}_{+}(t):\bw M_0\sra \left(\bw M_0\right)[[t]]$$ 
generated by $\sigma_{+1}$.  We shall denote by $\sigma_{+i},\ovsig_{+i}:\wM_0\sra \wM_0$ the coefficients of $t^i$ in $\sigma_+(t)$ and $\ovsig_+(t)$ respectively.
In the next subsection, the HS--derivations corresponding to $\sigma_{-1}$ will also appear.

Let us fix $r>0$. It is convenient to enumerate the basis of $\bw^rM_0$, which corresponds to $(b_j)_{j\geq 1}$,  by partitions $\blamb=(\lambda_1\geq\cdots\geq \lambda_r\geq 0)$ of lenght at most $r$. We write  $\Pcal_r$ for the set of all such partitions, and denote the basis vectors as follows,
\be
\wb^r_\blamb=b_{1+\lambda_r}\w b_{2+\lambda_{r-1}}\w\cdots\w b_{r+\lambda_1}\,, \ \ \blamb\in\Pcal_r\,.\label{brblamb}
\ee
In particular, the zero partition ${\bf 0}=(\lambda_1=0)$ gives $\wb^r_{\bf 0}=b_1\w \cdots \w b_r$.

\medskip
Now consider $(e_i)_{1\leq i\leq r}$ as indeterminates, and the polynomial ring 
\be
B_r={\QQ}[e_1,e_2,\ldots, e_r].\label{eq:Br}
\ee 
Let us equip $\bw^rM_0$ with a $B_r$-module structure via
\be 
E_r(t)\wb^r_\blamb=\overline{\sigma}_+(t)\wb^r_\blamb\,,\label{brmode}
\ee 
where $E_r(t)$ is given by (\ref{eqn:E}). In terms of the inverses, see~(\ref{eqn:H}), the same structure is given by
 \be 
H_r(t)\wb^r_\blamb=\sigma_+(t)\wb^r_\blamb\,.\label{brmodh}
\ee 

 The interpretation of  $e_j$'s and $h_j$'s as the elementary  and  the complete  symmetric functions of $r$ variables, see  Remark~\ref{Rem13} (2), suggests  to consider
 $B_r$ as a $\QQ$-vector spaces generated by the Schur polynomials  (see, for example,~\cite[I\,3]{MacDonald}),
\be
\Delta_\blamb(H_r)=\det(h_{\lambda_j-j+i})_{1\leq i,j\leq r}\,,\ \ \blamb\in\Pcal_r.\label{Schur}
\ee
Here $h_j$'s are defined as in Remark~\ref{Rem13} (1)  for $j\geq 0$, and $h_j=0$ for $j<0$.

According to {\sl Giambelli's formula} as in \cite[p.~321]{G1}), 
\be
\wb^r_\blamb=\Delta_\blamb(H_r)\wb^r_{\bf 0},\label{eq:giamb}
\ee
that is,  $\bw^rM_0$ is a free $B_r$-module of rank $1$ generated by $\wb^r_{\bf 0}$. This allows us  to equip $M_0$ with a multiplicative structure  over $B_r$, see Proposition~\ref{rModule}. 
Denote $M_0$, endowed with this  multiplicative structure, by $M_r$.  In  Section~\ref{S4}, we check that 
\begin{itemize}
\item  $M_r$ is a  $B_r$-module of rank $r$ freely generated by $b_1,\ldots, b_r$.

\item $\bw^rM_r$ is $\bw^rM_0$  with the $B_r$-module structure  defined by~(\ref{brmode}) or~(\ref{brmodh}). 

\item  $e_i$ is the eigenvalue of $\ovsig_{+i}$  restricted to $\bw^rM_r$, $1\leq i\leq r$.

\item  $h_j$ is the eigenvalue of the restriction of $\sigma_{+j}$ to $\bw^rM_r$, $j\geq 0$. 
\end{itemize}

\begin{remark}\label{RcharP} The notion of HS--derivation on an exterior algebra  enables one to extend some finite-dimensional linear algebra concepts (like eigenvalues and characteristic polynomials) to an  infinite-dimensional situation. Indeed, an endomorphism of an infinite-dimensional vector space does not have a characteristic polynomial, whereas  the corresponding HS--derivation is still  defined. 
\end{remark}
 
\subsection{Finite-dimensional approximations of bosonic vertex operators}\label{S23}
We apply the construction of the previous subsection in order to get 
 a ``finite-dimensional approximation''  of the well-known expression of the vertex operators occurring in the boson-fermion correspondence.  We interpret this approximation as the ratio of certain characteristic polynomials.
 
Details are in Section~\ref{S5}, see also~\cite{GatSal3}.

Take  the polynomial ring of countably many indeterminates, $B=\QQ[x_1,x_2,\ldots]$, and define the {\em bosonic vertex operators}, following~\cite[p.~56]{Ka},
\[
\Gamma(t)=\exp\left(\sum_{i\geq 1}x_it^i\right)\cdot\exp\left(-\sum_{i\geq 1}{1\over i t^i}{\d\over \d x_i}\right):B\sra B[t^{-1},t]],
\]
\[
\Gamma^*(t)=\exp\left(-\sum_{i\geq 1}x_it^i\right)\cdot\exp\left(\sum_{i\geq 1}{1\over i t^i}{\d\over \d x_i}\right):B\sra B[t^{-1},t]].
\] 

We find finite-dimensional counterparts of these operators using the symmetric functions interpretation. Namely, 
similarly to the finite-dimensional case, define $E_\infty(t)$ and $H_\infty(t)$,
\[E_\infty(t)=1-e_1t+e_2t^2+\cdots+(-1)^ke_kt^k+\ldots\ ,\ \ \ H_\infty(t)=1/E_\infty(t), \]
as the generating functions of the  elementary and the complete symmetric functions of a countable set of variables, say, $(\xi_k)_{k\geq 1}$.  Consider also $x_j=j\sum_{k\geq 1}\xi_k^j$,  {\sl the power sum  symmetric functions}, see~\cite[I\,3]{MacDonald}. We have $\QQ[x_1, x_2, \ldots ]=\QQ[e_1, e_2, \ldots ]$. Moreover,  $X_\infty(t)=\sum_{i\geq 1}x_it^i$,  the generating function of $(x_i)_{i\geq 1}$, satisfies
\[\exp\left(\sum_{i\geq 1}x_it^i\right)=\sum_{i\geq 0}h_it^i.\]
Clearly, $E_r(t)$, $H_r(t)$, $X_r(t)$ are obtained from $E_\infty(t)$,  $H_\infty(t)$,  $X_\infty(t)$ by setting $\tau_k=0$ for $k> r$. 

In order to define  $\Gamma_r(t)$ and $\Gamma^*_r(t)$ for $r>0$, use the notation of Section~\ref{S22}. In particular, the ring~(\ref{eq:Br}) is freely generated by  the Schur polynomials~(\ref{Schur}), {and $\bw^rM_r$ is spanned over $B_r$ by}  $\wb^r_{\bf 0}$, according to~(\ref{eq:giamb}).   Juxtaposing~(\ref{brmode}) or ~(\ref{brmodh}) and~(\ref{eq:giamb}) we get, respectively,
\be\overline{\sigma}_+(t)\wb^r_\blamb=\overline{\sigma}_+(t)\left(\Delta_\blamb(H_r)\wb^r_{\bf 0}\right)=E_r(t)\Delta_\blamb(H_r)\wb^r_{\bf 0}\label{eq:inducmap1},,\ee
\be\sigma_+(t)\wb^r_\blamb=\sigma_+(t)\left(\Delta_\blamb(H_r)\wb^r_{\bf 0}\right)=H_r(t)\Delta_\blamb(H_r)\wb^r_{\bf 0}\,.\label{eq:inducmap2}\ee

Thus each of $\overline{\sigma}_+(t)$, $\sigma_+(t)$  defines certain homomorphism $B_r\sra B_r[[t]]$,
which we denote in the same way. 

 For the HS--derivations generated by the shift operator  $\sigma_{-1}$ of Section~\ref{S22}, 
 we use the indeterminate $t^{-1}$ instead of $t$, and denote them by
$\sigma_-(t^{-1}),\ \overline{\sigma}_-(t^{-1})$. The corresponding homomorphisms are defined via 
\[
(\sigma_-(t^{-1})\Delta_\blamb(H_r))\wb^r_0=\sigma_-(t^{-1})\wb^r_\blamb,\quad (\overline{\sigma}_-(t^{-1})\Delta_\blamb(H_r)) \wb^r_0=\overline{\sigma}_-(t^{-1})\wb^r_\blamb.
\]
Definition of $\sigma_{-1}$ implies that 
$\sigma_-(t^{-1})b_i=b_i+b_{i-1}t^{-1}+b_{i-2}t^{-2}+\cdots+b_1t^{1-i}$
is a polynomial of $t^{-1}$  for each $i>0$. It follows  that  $\sigma_-(t^{-1}),\,\overline{\sigma}_-(t^{-1})$ also send  all $\Delta_\blamb(H_r)$'s to $B_r$-polynomials of $t^{-1}$. 

\medskip

Now we are ready to define the homomorphisms $\Gamma_r(t), \Gamma_r^*(t): B_r\sra B_r[t^{-1},t]]$ by their values on  $\Delta_\blamb(H_r)$'s, as follows,
\begin{eqnarray*}
\Gamma_r(t)\left(\Delta_\blamb(H_r)\right)&=&{1\over E_r(t)}\left(\overline{\sigma}_-(t^{-1})\Delta_\blamb(H_r)\right),\\
\Gamma^*_r(t)\left(\Delta_\blamb(H_r)\right)&=&E_r(t)\cdot\left(\sigma_-(t^{-1})\Delta_\blamb(H_r)\right).
\end{eqnarray*}
For $r_1<r_2$,  the natural projection  $B_{r_2}\to B_{r_1}$ sending each of $e_{r_1+1},\ldots, e_{r_2}$ to zero, sends $E_{r_2}(t)$ to $E_{r_1}(t)$, $H_{r_2}(t)$  to   
$H_{r_1}(t)$, and $X_{r_2}(t)$ to $X_{r_1}(t)$ . In this sense,  $E_r(t)\sra E_\infty(t)$, $H_r(t)\sra H_\infty(t)$, $X_r(t)\sra X_\infty(t)$ as $r\sra\infty$.

Thus, $\Gamma_r(t)$ and $\Gamma^*_r(t)$ tend to $\Gamma(t)$ and $\Gamma^*(t)$ when  $r\sra\infty$.

 \section{ Cayley--Hamilton Theorem revisited}\label{S3}
\subsection{Hasse-Schmidt derivations on  exterior algebras \cite{G1, SCGA}.} \label{Pre}  Let $A$ be  a commutative ring with unit,  $M$  a free $A$-module of rank $r$, and $b_1,\ldots,b_r$   some $A$-basis  of   $M$.

Set $\bw^0M=A$. For  $1\leq j\leq r$, denote by $\bw^jM$  the  $A$-module generated by all $b_{i_1}\wedge\ldots\wedge b_{i_j}$ modulo permutations,
\[
b_{i_{\tau(1)}}\w\ldots\w b_{i_{\tau(j)}}=\sgn(\tau)b_{i_1}\w\ldots\w b_{i_j},
\]
where $\sgn(\tau)$ is the sign of  permutation $\tau$.
In particular, $\bw^1M=M$. 

The exterior algebra $\wM=\bigoplus_{j=0}^r\bw^jM$ possesses the natural graduation given by  juxtaposition $\wedge:\bw^iM\times \bw^jM\sra \bw^{i+j}M$.

We  denote by $(\wM)[[t]]$ the ring of formal power series of $t$ with coefficients in $\wM$, and by $(\End_A(\wM))[[t]]$ the ring of formal power series of $t$ with coefficients in
$\End_A(\wM)$. 

For
$\Dcal(t)=\sum_{i\geq 0}\Dcal_it^i, \ 
{\Dcal}(t)=\sum_{j\geq 0}\widetilde{\Dcal}_jt^j\in(\End_A(\wM))[[t]]$,  their
product is defined as follows,
\[
\Dcal(t)\widetilde{\Dcal}(t)u=\Dcal(t)\sum_{j\geq 0}\widetilde{\Dcal}_ju\cdot t^j=\sum_{j\geq 0}(\Dcal(t)\tilde{\Dcal}_ju)\cdot t^j,\ \ \forall u\in\wM.
\]
Given series $\Dcal(t)$, we use the same notation for the induced  $A$-homomorphism,  
\[\Dcal(t):\wM\sra\wM[[t]],\ \  u\,\mapsto\, \Dcal(t)u=\sum_{i\geq 0}\Dcal_iu\cdot t^i,\ \ 
\forall u\in\wM.\]
The series $\Dcal(t)=\sum_{i\geq 0}\Dcal_it^i$   is {\em invertible} in  $(\End_A(\wM))[[t]]$, if there exists $\ovDc(t)\in(\End_A(\wM))[[t]]$ such that
\be \Dcal(t)\ovDc(t)=\ovDc(t)\Dcal(t)=\1_{\wM}.\label{eq:DovD}\ee
We call $\ovDc(t)$ the {\em  inverse} series and write it  in the form $\ovDc(t)=\sum_{i\geq 0}(-1)^i\ovDc_i t^i$. Then~(\ref{eq:DovD}) is equivalent to
\be
\Dcal_j-\ovDc_1\Dcal_{j-1}+\ldots+(-1)^j\ovDc_j=0, \ \forall j\geq 1. \label{eq:dinvdbar}
\ee
One can check that $\Dcal(t)$ invertible if and only if $\Dcal_0$ is an automorphism of $\wM$.

\begin{proposition}\label{equivl} {\em The following two statements are equivalent:} 
\begin{enumerate}
\item[$ i)$] $\Dcal(t)(u\w v)=\Dcal(t)u\w \Dcal(t)v$, $\ \forall u,v\in \wM$;
\item[$ii)$]  $\Dcal_i(u\w v)=\sum_{j=0}^i\Dcal_{j}u\w \Dcal_{i-j}v$, $\ \forall u,v\in \wM$, $\ \forall i\geq 0$.
\end{enumerate}
\end{proposition}
\proof
$i)\Rightarrow ii)$\ \ By definition of $\Dcal(t)$,  one can write $i)$ as
\be
\sum_{i\geq 0}\Dcal_i(u\w v)t^i=\sum_{j_1\geq 0}\Dcal_{j_1}u\cdot t^{j_1}\w\sum_{j_2\geq 0}\Dcal_{j_2}v\cdot t^{j_2}.\label{eq:CH02}
\ee
Hence $\Dcal_i(u\w v)$ is the coefficient of $t^i$ on the right hand side of~(\ref{eq:CH02}), which is  $\sum_{j_1+j_2=i}\Dcal_{j_1}u\w \Dcal_{j_2}v=\sum_{j=0}^i\Dcal_{j}u\w \Dcal_{i-j}v$.

\medskip

$ii)\Rightarrow i)$\ \  We have
\begin{eqnarray*}
\Dcal(t)(u\w v)&=&\sum_{i\geq 0}\Dcal_i(u\w v)t^i=\sum_{i\geq 0}\left(\sum_{i_1+i_2=j}\Dcal_{i_1}u\w \Dcal_{i_2}v\right)t^i\cr&=& \sum_{i\geq 0}\left(\sum_{i_1}\Dcal_{i_1}u\cdot t^{i_1}\w \sum_{i_2}\Dcal_{i_2}v\cdot t^{i_2}\right)=\Dcal(t)u\w \Dcal(t)v. \hskip90pt \qed
\end{eqnarray*}

\begin{definition} (Cf.~\cite{G1})  Let $\Dcal(t)\in (\End_A(\wM))[[t]]$. The induced map $\Dcal(t): \wM\sra(\wM)[[t]]$  is called a  {\em Hasse--Schmidt derivation} (or,  for brevity,  a {\em HS--derivation}) on  $\wM$, if it satisfies  the  (equivalent) conditions of Proposition~\ref{equivl}.
\end{definition}

\begin{proposition}\label{propHS} (Cf.~\cite{G1, {SCGA}})  {\em The product of two $HS$--derivations is a HS--derivation. The inverse of a HS--derivation is a HS--derivation.}
\end{proposition}
 \proof  For the product of HS--derivations  $\Dcal(t)$ and $\tilde{\Dcal}(t)$, the statement $i)$ of  Proposition~\ref{equivl} holds. Indeed, $\forall u,v\in\bw M$,
 \begin{eqnarray*}
\Dcal(t)\tilde{\Dcal}(t)(u\w v)\hskip-6pt&=&\hskip-6pt \Dcal(t)\left(\sum_{j\geq0}\ \sum_{j_1+j_2=j}\tilde{\Dcal}_{j_1}u\w\tilde{\Dcal}_{j_2}v\right)t^j\\
&=&\hskip-6pt\sum_{j\geq 0}\sum_{j_1+j_2=j}\Dcal(t) D_{j_1}u\cdot t^{j_1}\w \Dcal(t)\Dcal_{j_2}v\cdot t^{j_2}\\
&=&\hskip-6pt \Dcal(t)\tilde{\Dcal}(t)u\w \Dcal(t)\tilde{\Dcal}(t)v.\hskip66pt  
 \end{eqnarray*}
Similarly, if $\ovDc(t)$ is the inverse of the HS--derivation $\Dcal(t)$, then  $\forall u,v\in\bw M$,
 \begin{eqnarray*}
\hskip24pt \ovDc(t)(u\w v)&=&\ovDc(t)(\Dcal(t)\ovDc(t)u\w \Dcal(t)\ovDc(t)v)\cr
&=&(\ovDc(t)\Dcal(t))(\ovDc(t)u\w\ovDc(t)v)\\
&=&\ovDc(t)u\w\ovDc(t)v.\hskip 90pt\qed
 \end{eqnarray*}
 
\begin{corollary}\label{intbp} \cite{SCGA}\  If  $\ovDc(t)$ is the inverse of a HS--derivation $\Dcal(t)$, then 
\begin{eqnarray}
\Dcal(t)u\w v&=&\Dcal(t)u\w \Dcal(t)\ovDc(t)v=\Dcal(t)(u\w\ovDc(t)v),\cr
u\w \ovDc(t)v&=&\ovDc(t)\Dcal(t)u\w \ovDc(t)v=\ovDc(t)(\Dcal(t)u\w v)  \label{eq:intprt}
\end{eqnarray}
for all $u,v\in\bw M$.  Equivalenly, for any $k\geq 1$,
\begin{eqnarray}
\Dcal_ku\w v&=&\Dcal_k(u\w v)-\Dcal_{k-1}(u\w\ovDc_1v)+\ldots+(-1)^ku\w\ovDc_kv,\cr
u\w \ovD_kv&=&\ovD_k(u\w v)-\ovD_{k-1}(D_1u\w v)+\ldots+(-1)^kD_ku\w v\label{f25}
\end{eqnarray}
\end{corollary}

\subsection{Proof of the Theorem~\ref{thm23}}\label{s31} 
As we have seen in Proposition~\ref{prop11},
any endomorphism $f\in \End_A(M)$ defines two graded mutually inverse HS--derivations, 
$$\ovD(t)=\bw\left(\1-ft\right)\ \ {\rm and}\ \  D(t)=\bw\left(\sum_{i\geq 0}f^it^i\right),$$ 
where $\1$ denotes the identity endomorphism.
Write  $\ovD(t)$ and $D(t)$ in the form 
$$\ovD(t)=\sum_{i\geq 0}(-1)^i\ovD_it^i \ \ {\rm and}\ \  D(t)=\sum_{i\geq 0}D_it^i,$$
then  these HS--derivations satisfy the following properties.
\begin{lemma}\label{obs} We have 

$(i)$  $\ovD_0=D_0=\1_{\bw M}$ and $D_1=\ovD_1$,

$(ii)$  ${D_i}\mid_M=f^i$, $i\geq 0$, 

$(iii)$ $\ovD_ku=0$, for all $u\in \bw^iM$ with $i<k$.
\end{lemma}
Indeed, $\ovD(t)\mid_{\bw^k M}$ is a polynomial of $t$ of degree $k$,  $1\leq k\leq r$. \qed

As before, we assume that our $A$-module $M$ is freely generated by $(b_j)_{1\leq j\leq r}$.   
Thus $\bw^rM$ has rank 1 and is spanned by $\wb^r_{\bf 0}=b_1\w \cdots \w b_r$. 
The restriction of each
$\ovD_i$ to  $\bw^rM$ is a multiplication by some scalar $e_i\in A$,
\be \ovD_i\left(\wb^r_{\bf 0}\right)=e_i\wb^r_{\bf 0},\ \ 1\leq i\leq r.\label{eibr}\ee
Take now $u\in\bw^iM$ and $ v\in\bw^{r-i}M$. Then  $D_ju\w v\in \bw^rM$ for $1\leq j\leq k$. 
Applying~(\ref{f25}) to our situation, we can write
\[D_ku\w v-e_1(D_{k-1}u\w v)+\cdots+(-1)^ke_k(u\w v)=(-1)^ku\w \ovD_kv\]
for $1\leq k\leq r$, and
\[D_ku\w v-e_1(D_{k-1}u\w v)+\cdots+(-1)^ke_r(D_{k-r}u\w v)=(-1)^ku\w \ovD_kv
\]
for $k>r$.

Equivalently, we have 
\be\left(D_ku-e_1D_{k-1}u+\cdots+(-1)^ke_ku\right)\w v=(-1)^ku\w \ovD_kv,\ 1\leq k\leq r, \label{eq:3}\ee
and
\be\left(D_ku-e_1D_{k-1}u+\cdots+(-1)^ke_rD_{k-r}u\right)\w v=(-1)^ku\w \ovD_kv,\ k>r. \label{eq:4}\ee
Of course, one can set $e_k=0$ for $ k>r$, in order do not distinguish between the two cases. However, we prefer a division into cases.

Assume now $i>r-k>0$. Then, according to Remark~\ref{obs}, $(iii)$, the right hand side of~(\ref{eq:3}) vanishes $\forall v\in\bw^{r-i}M$, as $r-i<k$.  This means that 
$$
D_ku-e_1D_{k-1}u+\cdots+(-1)e_ku=0
$$
for any $u\in\bw^iM$ with $i>r-k>0$. This proves the first part of Theorem~\ref{thm23}. 

If $k>r$, then the left hand side of~(\ref{eq:4}) vanishes for each $i\geq 0$, and this proves the second part. \hskip3truecm \qed 

\begin{remark}
Thus we understand (\ref{eq:1}) as an abstract Cayley--Hamilton theorem valid for general invertible HS--derivations on exterior algebras of arbitrary  (not necessarily free) modules. If the module is free and at most countably generated, then~(\ref{eq:1}) produces a sequence of Cayley--Hamilton relations~(\ref{f25}). This sequence turns into the classical Cayley--Hamilton formulas~(\ref{eqt22:prech}) when the HS--derivation corresponds to an endomorphism of a finitely generated free module.
\end{remark}

\subsection{Example} Take $M=\RR^3$.  Let $f:\RR^3\sra\RR^3$ have an eigenbasis, 
\[f\bfu = a\bfu,\ \ f\bfv=b\bfv, \ \ f\bfw=c\bfw.\]
Then $\det(\1t-f)=t^3-e_1t^2+e_2t-e_3$, where
\[e_1=a+b+c,\ \ e_2=ab+ac+bc, \ \ e_3=abc.
\]
In notation of Theorem~\ref{thm23}, we have $r=3$.

$1)$\ \ Let us take $k=2$ and check that $D_2-e_1D_1+e_2\1$ vanishes on  $\RR^3\w\RR^3$, 
as it should be, according to~(\ref{eq2:prech}). 
We show the calculation of $(D_2-e_1D_1+e_2\1)(\bfu\w \bfv)$; for two other basis vectors $\bfu\w \bfw$ and $\bfv\w \bfw$ it is completely similar.
 
First, we find the action of $D_1, D_2$ on $\bfu\w \bfv$. We have
\[(\1+ft+f^2t^2+\circ(t^2))\bfu\w (\1+ft+f^2t^2+\circ(t^2))\bfv=\]
\[=\bfu\w \bfv+(f\bfu\w \bfv+\bfu\w f\bfv)t+(f^2\bfu\w \bfv+f\bfu\w f\bfv+\bfu\w f^2\bfv)t^2+\circ(t^2).
\]
Thus\ \  $D_1(\bfu\w \bfv)=(a+b)(\bfu\w \bfv)$,\ \  $D_2(\bfu\w \bfv)=(a^2+ab+b^2)(\bfu\w \bfv)$,
and 
\[(D_2-e_1D_1+e_2\1)(\bfu\w \bfv)=(a^2+ab+b^2-e_1a-e_1b+e_2)(\bfu\w \bfv).\]

Now, we substitude the expressions for $e_1$, $e_2$, and get
\[a^2+ab+b^2-e_1a-e_1b+e_2=a^2+ab+b^2-a^2-ab-ac-ab-b^2-bc+ab+ac+bc=0.\]

$2)$\ \ Take $k=4$ and check that $D_4-e_1D_3+e_2D_2-e_3D_1$ vanishes on $\RR^3\w\RR^3$.
According to~(\ref{eqt22:prech}), this endomorphism vanishes on the whole of $\bw\RR^3$,
and in fact the verification for the rest of the direct summands is simpler. Again we calculate the image of $\bfu\w \bfv$.

First, we obtain $D_3(\bfu\w \bfv)$ and $D_4(\bfu\w \bfv)$ in the standard way, writing
\[(\1+ft+f^2t^2+f^3t^3+f^4t^4+\circ(t^5))\bfu\w (\1+ft+f^2t^2+f^3t^3+f^4t^4+\circ(t^5))\bfv,\]
and collecting the coefficients of $t^3$, $t^4$, respectively. We get
\[D_3(\bfu\w \bfv)=(a^3+a^2b+ab^2+b^3)(\bfu\w \bfv), D_4(\bfu\w \bfv)=(a^4+a^3b+a^2b^2+ab^3+b^4)(\bfu\w \bfv),
\]
substitute all the expressions for $D_4, D_3, D_2, D_1, e_1,\ e_2,\ e_3$ in terms of $a,b,c$ into $D_4-e_1D_3+e_2D_2-e_3D_1$, and safely get $0$.

\begin{remark} \ \ In general, if $f\in\End(M)$ is diagonalizable, and if  $(\bfv_i)_{1\leq i\leq r}$ is an eigenbasis,  $f\bfv_i=x_i\bfv_i$, $1\leq i \leq r$, then the vector $\bfv_1\w\ldots\w\bfv_l\in\bw^lM$ is an eigenvector of $D_k$ with the eigenvalue which is the complete symmetric polynomial of $x_1,\ldots, x_l$ of degree $k$. 
Therefore, for a diagonizable endomorphism our Theorem~\ref{thm23} is reduced to the following identity. 
Denote by $h_i(\xx_j)$ the complete symmetric polynomial of degree $i$ in $x_1,\ldots,x_j$, and by $e_k(\xx_n)$ the elementary symmetric polynomial of degree $k$ in $x_1,\ldots,x_n$. Then, for $n\geq 1$ and all $1\leq j\leq n$, we have
\[h_n(\xx_j)-e_1(\xx_n)h_{n-1}(\xx_j)+\ldots+(-1)^ne_n(\xx_n)=0.
\]
One can deduce the identity, for example,  from the formula (*) of~\cite[I 3 28]{MacDonald}. 
\end{remark}

This remark can be turned into a rigorous general proof, using a standard (though rather long) reasoning. Another possible way, which was suggested by our referee, is based on the Frobenius proof of the classical Cayley--Hamilton theorem for the complex matrices, \cite{Fro}. We were not aware of that 1896 paper by Frobenius. Probably one could translate our arguments into the language of matrix minors. 
However, our approach, through the relationship to symmetric functions, is short, easy, and, in addition, allows us to concern with the infinite-dimensional case.

\subsection{The case of a $\QQ$-algebra.}\label{s34}  If $A$ is a $\QQ$-algebra, then for $f\in \End_AM$ the exponential 
\[
\exp(ft):=\sum_{k\geq 0}{{f^kt^k}\over{k!}}\in (\End_AM)[[t]]
\]
is well-defined.
In the ring $(\End_AM)[[t]]$, there is the formal derivative by $t$,
\[ y(t)=\sum_{k\geq 0} g_kt^k \ \Rightarrow\  y'(t)=\sum_{k\geq 0} kg_kt^{k-1},\ \ g_k\in\End_AM.
\]
Recall the notation  $E_r(t)$ given by~(\ref{eqn:E}),  and write the characterictic polynomial of $f$ as in~(\ref{eqn:1E}).
In~\cite{gatlak} for any commutative ring $R$ containing the rational numbers, {\em the formal Laplace transform} $L:R[[t]]\sra R[[t]]$  and its inverse $L^{-1}$ are defined as follows, 
\[
L\,\sum_{n\geq 0}a_nt^n=\sum_{n\geq 0}n!a_nt^n,\ \ \
L^{-1}\sum_{n\geq 0}c_nt^n=\sum_{n\geq 0}c_n{t^n\over n!},\ \ a_n, c_n\in R.
\] 

Take the  {\sl inverse formal Laplace transform  of the HS--derivation} $D(t)=\sum_{i\geq 0}D_it^i$ corresponding to
$f\in \End_AM$,
\be
D^*(t)=L^{-1}\sum_{k\geq 0}D_kt^k=\sum_{k\geq 0}{{D_kt^k}\over{k!}}\in \left(\End_A(\bw M)\right)[[t]].\label{D*}
\ee
Define $\ttp_k(D)$ as the coefficient of $t^k$ in $ E_r(t)D(t)$. We have
\[ \ttp_0(D)=\1, \ \ttp_1=D_1-e_1\1,\ \ttp_j=D_j-e_1D_{j-1}+\ldots+(-1)^je_j\1,\ 1<j<r,\]
and
\[\ttp_{r+j}(D)=D_{r+j}-e_1D_{r+j-1}+\ldots+(-1)^re_rD_j=\0,\ \ j\geq 0,\]
according to Theorem~(\ref{thm23}). Therefore,

\begin{proposition} We have
\be
D(t)={\1+\ttp_1(D)t+\ldots+\ttp_{r-1}(D)t^{r-1}\over E_r(t)}. \qquad \ \ \qquad \qed \label{eq:prexp}
\ee
\end{proposition}

\begin{corollary}\label{coru} Let $\QQ\subseteq A$ and the characteristic polynomial of $f\in\End_AM$ be given by~(\ref{eqn:1E}). Then the series $D^*(t)$ defined in~(\ref{D*}) solves
the ordinary differential equation 
\be
y^{(r)}(t)-e_1y^{(r-1)}(t)+\ldots+(-1)^re_ry(t)=0\label{ODE}  
\ee 
 in $(\End_A(\wM))[[t]]$.
\end{corollary}
\proof  Take the inverse formal Laplace transform of~(\ref{eq:prexp}). We obtain
\[
D^*(t)=u_0+\ttp_1(D)u_{-1}+\ldots+\ttp_{r-1}(D)u_{-r+1},
\]
where
\[
u_{-j}=u_{-j}(t)=L^{-1}\left(t^j\over E_r(t)\right), \ \ 0\leq j\leq r-1.
\]
Let us re-write the series $u_0,u_{-1},\ldots,u_{-r+1}$ in terms of $H_r(t)=1/E_r(t)$, see Remark~\ref{Rem13}(1), 
\[
u_{-j}=L^{-1}(t^jH_r(t))=\sum_{n\geq j}h_{n-j}{t^n\over n!},\ \ 0\leq j\leq r-1.
\]
In~\cite{GaSch1}, we proved that these series form an $A$-basis of solutions to the  ODE (\ref{ODE}) in $R[[t]]$. For $R=\End_A(\wM)$ we get the claim.  \qed

\subsection{Elementary remarks.}
We finish this section with a few remarks relevant to the case when  $A$ is a $\QQ$-algebra.

(1)\ \ The characteristic polynomial of $f\in\End_AM$ is given by~(\ref{eqn:1E})  if and only if $y(t)=\exp(ft)$ satisfies the linear ordinary differential equation (\ref{ODE}). This is our  Corollary~\ref{coru} restricted to $M$.

 In particular,  
 \[\exp(ft)=v_0(t)\1_M+v_1(t)f+\cdots+v_{r-1}(t)f^{r-1},\]
where  $(v_j(t))_{0\leq j\leq r-1}$ is  the {\em standard} $A$-basis of solutions  to (\ref{ODE}) in $A[[t]]$ , that is,  $v_j^{(i)}(t)=\delta_{ij}$, $0\leq i,j\leq r-1$. 
Indeed, $\1_M, f, \ldots, f^{r-1}$ are the initial conditions of the solution $\exp(ft)$.

In the context of endomorphisms of complex vector spaces, the formula for $\exp(ft)$ was obtained  in 1966 by Putzer \cite{Pu}, and then re-obtained in 1998 by  Leonard and Liz, \cite{Leon1, Liz}, in a different way.

(2)\ \ The relation between the standard fundamental system $v_j(t))_{0\leq j\leq r-1}$ and the fundamental system  ${u_{-j}(t)}_{0\leq j\leq r-1}$ appeared in the proof of Corollary~\ref{coru} is
as follows. Consider the linear system of first order differential equations equivalent to our ODE (\ref{ODE}),
\[ y_1'=y_2,\ \ y_2'=y_3,\  \ldots,\  y_{r-2}'=y_{r-1},\ \ y_{r-1}'=e_1y_{(r-1)}-\ldots+(-1)^{r-1}e_ry_1. \]
Denote the matrix of this system by $P_r$. Then $Q=\exp(P_rt)$ is the Wronski matrix of  $v_1(t), \ldots, \ v_{r-1}(t)$,
\[ (Q)_{ij}= v_j^{(i)}(t), \ \ 0\leq i,j\leq r-1,\]
and $u_0(t), u_{-1}(t),\ldots, \ u_{1-r}(t)$ is the last column of $Q$,
\[v_{r-1}(t)=u_{1-r}(t),\ v_{r-1}'(t)=u_{2-r}(t), \ldots,\ 
v_{r-1}^{(r-1)}(t)=u_0(t).\]

(3)\ \ As another elementary corollary of our considerations, we get formulas for the coefficients $e_k$ of the characterictic polynomial of $f\in\End_AM$ in terms of its matrix elements. If $C=(c_{ij})$ is the $r\times r$ matrix of $f$ in some $A$-basis of $M$, denote by $D(i_1,\ldots,i_k)$ the determinant of the $(r-k)\times (r-k)$-matrix obtained from the matrix $C$ by deleting the $i_1$-th,\ldots,$i_k$-th rows and columns.  Then 
\[(-1)^ke_k=\sum_{1\leq i_1<\ldots <i_k\leq r}D(i_1,\ldots,i_k),\ 1\leq k\leq r.
\]
This formula was obtained differently by Brooks in~\cite{Brooks}.

\section{Countably generated $\QQ$-vector spaces}\label{S4} 
Let $M_0$ be a  $\QQ$-vector space generated by $(b_j)_{j\geq 1}$ and $\bw M_0=\bigoplus_{r\geq 0}\bw^rM_0$ be its exterior algebra. 

As in~Section \ref{S22}, we take shift operators $\sigma_{\pm 1}\in \End_\QQ M_0$, and
denote by $\sigma_{\pm}(t^{\pm 1}), \overline{\sigma}_{\pm}(t^{\pm 1})$ the corresponding $HS$-derivations. Below we will write  
\be\sigma_+(t)=\sum_{j\geq 0}\sigma_jt^j, \ \ \ 
\overline{\sigma}_+(t)=\sum_{j\geq 0}\overline{\sigma}_jt^j,\label{sigma}\ee
skipping sign $+$ in the subscript.

In this section, we treat  $e_1,\ldots, e_r$ as indeterminates. As we already pointed out in Section~\ref{S22}, the ring $B_r$, given by (\ref{eq:Br}), has a basis formed by Schur polynomials $\Delta_\blamb(H_r)$, see~(\ref{Schur}). The  structure of a principal $B_r$-module on
$\bw^rM_0$ is defined via any of the two equivalent equalities~(\ref{brmode}) and~(\ref{brmodh}).

Let $(\beta_i)_{i\geq 1}$ be linear forms on $M_0$ defined by $\beta_i(b_j)=\delta_{ij}$. Their linear span is, by definition, the restricted dual $M_0^*$. Each $\beta_j$ induces a $\QQ$-linear contraction map $\beta_j:\bw^rM_0\sra \bw^{r-1}M_0$ defined by $\beta_j\lrcorner m=\beta_j(m)$ for all $m\in M_0$ and
\be 
\beta_j\lrcorner (m\w \eta)=\beta_j(m)\eta-m\w \beta_j\lrcorner \eta.\label{eq:basisc}
\ee 
As each $\zeta\in \bw M_0$ is a sum of homogeneous elements of the form $m\w \eta$, equation~(\ref{eq:basisc}) defines the contraction operator over the entire exterior algebra $\bw M_0$.

\begin{lemma} Let $m,m'\in M_0$ satisfy 
\be
m\w\eta=m'\w\eta\label{eq:hyp}
\ee for all $\eta\in \bw^{r-1}M_0$. Then $m=m'$.
\end{lemma}
\proof 
Under the hypothesis (\ref{eq:hyp}), suppose first that $m'=am$ for some $a\neq 1$. If $m\neq 0$, then there is $\mu\in M_0^*$ such that $\mu(m)\neq 0$. Because of the isomorphism $\bw^{r-1}M_0\cong M_0^*$, there is then $\eta\in \bw^{r-1}M_0$ such that $m\w\eta\neq 0$. We get $m'\w \eta=a(m\w\eta)$,  hence $m\w\eta\neq m'\w\eta$.   

If $m$ and $m'$ are not proportional, take their duals $\mu,\mu'\in M_0^*$  and choose $\zeta_{m,m'}=\mu'\lrcorner (\mu\lrcorner\wb^r_0)\in \bw^{r-1}M_0$.  Then  $\mu(m)=\mu'(m')=1$ and so $m\w m\w\zeta_{m,m'}=0$ while $m'\w m\w \zeta_{m,m'}\neq 0$.\qed

\begin{corollary}\label{em} 
For each $e_i$,  $1\leq i\leq r$, there exists the unique mapping $M_0\sra M_0$, called the {\sl multiplication by $e_i$},  $m\mapsto e_im$, such that
\be
(e_im)\w \eta=e_i(m\w \eta),\label{eq1:Mr}
\ee
for all $\eta\in \bw^{r-1}M_0$.
\end{corollary}
The vector space $M_0$ endowed with the multiplications by $e_i$, $1\leq i\leq r$, becomes a $B_r$-module, denoted by $M_r$.   
 
\begin{proposition}\label{rModule}
The  $B_r$-module $M_r$ is freely generated by $(b_i)_{1\leq i\leq r}$, and $\sigma_{+}(t)\in End_{B_r}(M_r)[[t]]$.
\end{proposition}
\proof
Let us check first that $b_1,b_2,\ldots, b_{r}$ are $B_r$-linearly independent. 
Denote by  $(\beta_j)_{j\geq 0}$ the generators of $M_0^*$ dual to $(b_i)_{i\geq 1}$.
Notice that
$a_1b_1+\cdots+a_{r}b_{r}=0$ implies $0=a_ib_i\w \eta_i\,,$
where $\eta_i=\beta_i\lrcorner \wb^r_{\bf 0}$, that is, $a_i=0$  for all $1\leq i\leq r$.

Now, let us show that $b_{i+r}-e_1b_{i+r-1}+\cdots+e_rb_{i}=0$ for all $i\geq 0$. 
This will prove, by induction, that $M_r$ is generated over $B_r$ by $b_1,b_2,\ldots, b_{r}$.  It is enough to observe that 
\[
\sum_{i=1}^\infty (b_i-e_1b_{i-1}+(-1)^re_rb_{r-1+i})t^i=E_r(t)\sigma_+(t)b_1
\]
is a polynomial of degree $r$ (here we set $b_j=0$ for $j<1$).  By definition of the module structure, for each $\eta \in \bw^{r-1}M_0$ we have
\begin{center}
\begin{tabular}{rrlr}
$E_r(t)(\sigma_+(t)b_0)\w \eta$&$=$&$E_r(t)(\sigma_+(t)b_0\w \eta)$& (definition of the\\
&&& $B_r$-module structure)\\ \\
&$=$&$E_r(t)\sigma_+(t)(b_0\w \ovsig_+(t)\eta)$&(Proposition~\ref{intbp})\\ \\
&$=$&$E_r(t)\displaystyle {1\over E_r(t)}(b_0\w \ovsig_+(t)\eta)$&(definition of the\\
&&& $B_r$-module structure)\\ \\
&$=$&$b_0\w \ovsig_+(t)\eta.$
\end{tabular}
\end{center}
We use now the agreement (\ref{sigma}). As in Remark~\ref{obs}, we see that $\ovsig_r$ vanishes on $\bw^{r-1}M_0$, hence the expression obtained above is a polynomial in $t$ of degree $r-1$.
We have so proven that $M_r$ is a $B_r$-module of rank $r$. Moreover $\sigma_1$ is 
$B_r$-linear.  In fact,
\[
\sigma_1(e_im)\w \eta=\sigma_1(e_im\w \eta)-e_im\w \sigma_1\eta=e_i\sigma_1(m\w \eta)-e_im\w\sigma_1\eta=
e_1\sigma_1m\w \eta. \qquad \qed
\]

\begin{corollary} The elements $e_i$, $1\leq i\leq r$,  and $h_j$, $j\geq 0$, of $B_r$ are  the eigenvalues of $\overline{\sigma}_i$ and  ${\sigma_j}$, respectively, thought of as endomorphisms of $\bw^rM_r$. \qed
\end{corollary}
Similarly to (\ref{eq:inducmap1}) and (\ref{eq:inducmap2}), the mappings $\sigma_-(t^{-1})$ and $\overline{\sigma}_-(t^{-1})$ define two   homomorphisms $B_r\sra B_r[t^{-1}]$, via the equalities
\be
(\sigma_-(t^{-1})\Delta_\blamb(H_r))\wb^r_0:=\sigma_-(t^{-1})\wb^r_\blamb,\label{feq:sigm_-}
\ee
and
\be
(\overline{\sigma}_-(t^{-1})\Delta_\blamb(H_r)) \wb^r_0:=\overline{\sigma}_-(t^{-1})\wb^r_\blamb.\label{feq:ovsigm_-}
\ee
By  abuse of notation, we denote the homomorphisms in the same way.

\section{Bosonic Vertex Operators}\label{S5}
\subsection{}
Let $B:=\QQ[x_1,x_2,\ldots]$ be the polynomial ring in infinitely many indeterminates and $M_0:=\bigoplus_{i\geq 0}\QQ b_i$.
 The aim of this section is to show that the {\em bosonic vertex operators}
\[
\Gamma(t):=\exp(\sum_{i\geq 1}x_it^i)\cdot\exp\left(-\sum_{i\geq 1}{1\over i t^i}{\d\over \d x_i}\right):B\sra B[t^{-1},z]]
\]
and 
\[
\Gamma^*(t):=\exp(-\sum_{i\geq 1}x_it^i)\cdot\exp\left(\sum_{i\geq 1}{1\over i t^i}{\d\over \d x_i}\right):B\sra B[t^{-1},z]]
\] 
may be identified with  ratios of characteristic series operators associated to the shift endomorphisms of step $\pm 1$ of $M_0$.

Let
$
\Gamma_r(t), \Gamma_r^*(t): B_r\sra B_r[t^{-1},t]]
$
be defined by
\[
\Gamma_r(t)\Delta_\blamb(H_r)\wb^r_0:=\sigma_+(t)\overline{\sigma}_-(t^{-1})\Delta_\blamb(H_r))\wb^r_0
\]
and
\[
\Gamma^*_r(t)\Delta_\blamb(H_r)\wb^r_0:=\overline{\sigma}_+(t)\sigma_-(t^{-1})\Delta_\blamb(H_r)\wb^r_0.
\]
Then, due to~(\ref{eq:inducmap1}) and~(\ref{eq:inducmap2}) one can write:
\[
\Gamma_r(t)={1\over E_r(t)}\cdot \overline{\sigma}_-(t^{-1})\qquad\mathrm{and}\qquad \Gamma_r^*(t)=E_r(t)\cdot {1\over\overline{\sigma}_-(t^{-1})}\,.
\]
\begin{remark} Notice that   $E_r(t)$ is indeed the characteristic polynomial of $\sigma_1$, thought of  as endomorphism of $M_r$,  and $\ovsig_-(t^{-1})$ is the characteristic series operator associated to $\sigma_{-1}$.
\end{remark}
\begin{proposition} \cite{GatSal3}\ \
The operators $\Gamma_r(t),\Gamma^*_r(t)$ tend to $\Gamma(t), \Gamma^*(t)$ as $r$ goes to infinity. 
\end{proposition}
\proof We sketch the arguments of  \cite{GatSal3}.
First of all, notice that for all $r\geq 1$
\be
\ovsig_-(t^{-1})h_n=h_{n}-{h_{n-1}\over {t}}\quad \mathrm{and}\quad \sigma_-(t^{-1})h_n=\sum_{i\geq 0}{h_{n-i}\over t^i}\,.\label{eq3:vohn}
\ee
Let  $\ovsig_-(t^{-1})H_r=(\sigma_-(t^{-1})h_n)_{n\in\ZZ}$ and $\sigma_-(t^{-1})H_r=(\sigma_-(t^{-1})h_n)_{n\in\ZZ}$. Then 
\be
\ovsig_-(t^{-1})\Delta_\blamb(H_r)=\Delta_\blamb(\ovsig_-(t^{-1})H_r), \ \  \sigma_-(t^{-1})\Delta_\blamb(H_r)=\Delta_\blamb(\sigma_-(t^{-1})H_r),\label{eq3:vodhn}
\ee
by~\cite[Propositions 6.2.10 and  6.2.13]{GatSal3}. Now, according 
to~\cite[Corollaries 6.2.11 and 6.2.14]{GatSal3},
\be
\sigma_-(t^{-1})(h_{i_1}\cdot\ldots\cdot h_{i_r})=\sigma_-(t^{-1})h_{i_1}\cdot\ldots\cdot\sigma_-(t^{-1}) h_{i_r}\label{eq:prehom1}
\ee
and
\be
\ovsig_-(t^{-1})(h_{i_1}\cdot\ldots\cdot h_{i_r})=\ovsig_-(t^{-1})h_{i_1}\cdot\ldots\cdot\ovsig_-(t^{-1}) h_{i_r}. \label{eq:prehom2}
\ee
Clearly formulas~(\ref{eq3:vohn}) and~(\ref{eq3:vodhn}) do not depend on $r$ when $r$ is big enough (that is, at least the length of the partition $\blamb$). Thus these formulas hold for $r=\infty$ as well.  We set 
\[E_\infty(t)=1-e_1z+e_2t^2+\cdots\quad \mathrm{and}\quad H_\infty(t)=1/E_\infty(t).\]  
Define now 
\[\exp\left(\sum_{i\geq 1}x_it^i\right):=\sum_{n\geq 0}h_nt^n.\]
Then we have
\[
B=\QQ[e_1,e_2,\ldots]=\QQ[h_1,h_2,\ldots]=\QQ[x_1,x_2,\ldots],
\]
see, for example,~\cite[I\,3]{MacDonald}.
Moreover,  
\be
{\partial h_n\over \partial x_i}=h_{n-i}\quad \mathrm{and}\quad
{\d^j h_n\over \d x_1^j}={\d h_n\over \d x_j}.\label{eq:dx1dxi}
\ee
As it follows from~(\ref{eq:prehom1}) and ~(\ref{eq:prehom2}), $\ovsig_-(t^{-1})$ and $\sigma_-(t^{-1})$,  for $r=\infty$, become ring homomorphisms $B\sra B[t^{-1}]$.
Thus
\begin{center}
\begin{tabular}{rlll}
$\ovsig_-(t^{-1})h_n$&$=$&$h_n-\displaystyle{h_{n-1}\over {t}}$&\hskip20pt  (first formula in~(\ref{eq3:vohn}))\\ \\
&$=$&$\left(1-\displaystyle{1\over t}\displaystyle{\partial\over\partial x_1}\right)h_n$ \\ \\
&$=$&$\exp\left(\displaystyle{-\sum_{i\geq 1}{1\over it}{\d^i\over \d x_1^i}}\right)h_n$&\hskip20pt  (definition of the logarithm \\
&&&\hskip22pt   of a formal power series)\\ \\
&$=$&$\exp\left(\displaystyle{-\sum_{i\geq 1}{1\over it}{\d\over \d x_i}}\right)h_n$ &\hskip20pt  (second equality in~(\ref{eq:dx1dxi})).
\end{tabular}
\end{center}

Notice that $\exp\left(\displaystyle{-\sum_{i\geq 1}{1\over it}{\d\over \d x_i}}\right)$, being  the exponent of a first order differential operator, is a ring homomorphism whose value at $h_n$ coincides with $\ovsig_-(t^{-1})h_n$.  This means that
\[
\ovsig_-(t^{-1})=\exp\left(\displaystyle{-\sum_{i\geq 1}{1\over it}{\d\over \d x_i}}\right).
\]
Similarly one shows that
\[
\sigma_-(t^{-1})=\exp\left(\displaystyle{\sum_{i\geq 1}{1\over it}{\d\over \d x_i}}\right)
\]
Thus
\[
\Gamma_\infty(t)={1\over E_\infty(t)}\ovsig_-(t)=\exp\left(\sum_{i\geq 1}x_it^i\right)\exp\left(\displaystyle{-\sum_{i\geq 1}{1\over it}{\d\over \d x_i}}\right)=\Gamma(t)
\]
and
\[
\Gamma^*_\infty(t)={E_\infty(t)}\sigma_-(t^{-1})=\exp\left(-\sum_{i\geq 1}x_it^i\right)\exp\left(\displaystyle{\sum_{i\geq 1}{1\over it}{\d\over \d x_i}}\right)=\Gamma^*(t)
\]
as claimed. \qed

\bibliographystyle{amsalpha}

\end{document}